\documentclass[11pt]{article}
\usepackage{amsmath, amsthm, amssymb}   
\usepackage{bridges}		
\usepackage{graphicx}			
\usepackage[colorlinks=true, urlcolor=blue, citecolor=black, linkcolor=black]{hyperref}  
\usepackage{subcaption}	
\usepackage{wrapfig}
\usepackage{tikz}
\usepackage{float}

\DeclareMathOperator\link{link}

\newtheorem*{definition*}{Definition}

\title{Knitting Knots \& the Framing Anomaly}

\author{Nadav Drukker\textsuperscript{1}, Elise Paznokas\textsuperscript{2}, and Dominik Schrimpel\textsuperscript{3}
\vspace{10pt}\\
\textsuperscript{1}Department of Mathematics, King’s College London, London, WC2R 2LS, United Kingdom; \\ \href{mailto:nadav.drukker@gmail.com}{nadav.drukker@gmail.com}\\
\textsuperscript{2}School of Mathematics, Trinity College Dublin, College Green, Dublin 2, Ireland \\ \href{mailto:epazno@gmail.com}{epazno@gmail.com} \\
\textsuperscript{3}École Polytechnique, Rte de Saclay, 91120 Palaiseau, France \\ \href{mailto:dominik.schrimpel@kcl.ac.uk}{dominik.schrimpel@kcl.ac.uk}}

\date{}					

\begin{document}

\maketitle

\thispagestyle{empty}

\vspace{-1cm}
\begin{abstract}
We study the twisting fault emerging in circular knitting, and its relation to the mathematical 
concepts of framing curves and the Gauss linking integral. We create three knitted bands with 
framing zero, one, and negative two, and use three different techniques to compute the framing 
using the Gauss linking integral. We also briefly mention the connection to in Chern–Simons gauge theory.
\end{abstract}

\vspace{-0.5cm}
\section*{Introduction}
One should not be surprised that there is a rich connection between the ancient craft of knitting and the mathematics of knot theory. See for example \cite{Belcastro, Seaton, Matsumoto}. 
One of the most elegant formulations of knot theory relies of the theoretical physics notion of topological quantum field theory, pioneered by Witten \cite{witten}. 
In this work we point out a less obvious relation between knitting and Witten's work: the 
twisting fault arising in circular knitting and the concept of framing arising when computing 
the expectation value of the Wilson loop observable in Chern-Simons gauge theory. While this 
is the motivation for this paper, we focus mostly on the relation between knitting, framing 
and the Gauss linking number.

\begin{figure}[H]
    \centering
    \begin{minipage}[b]{0.49\textwidth}
	    \includegraphics[width=\textwidth]{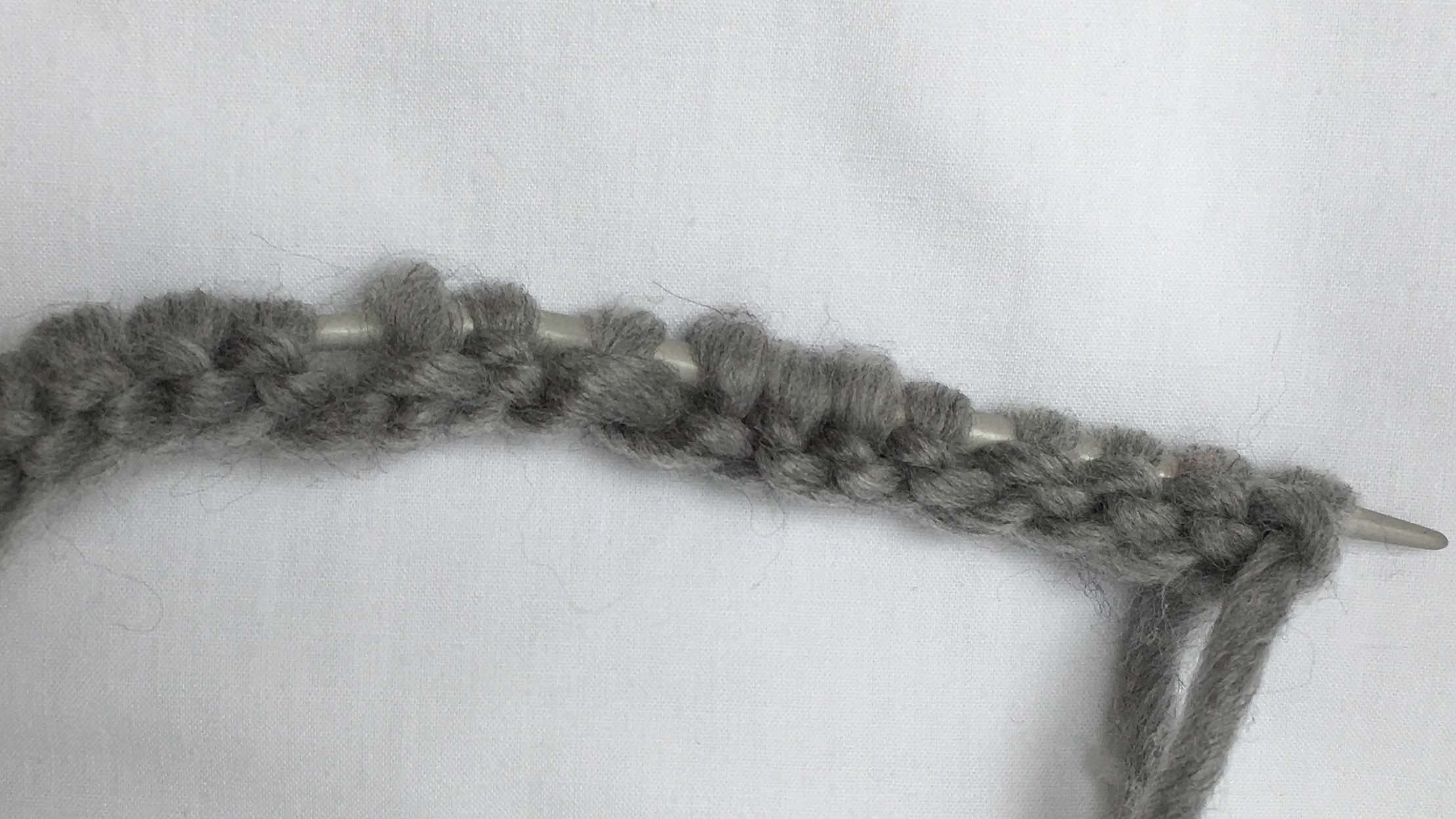}
       	    \subcaption{circular needles with cast-on}
        	\label{fig:1a}
    \end{minipage}
~ 
    \begin{minipage}[b]{0.49\textwidth} 
	    \includegraphics[width=\textwidth]{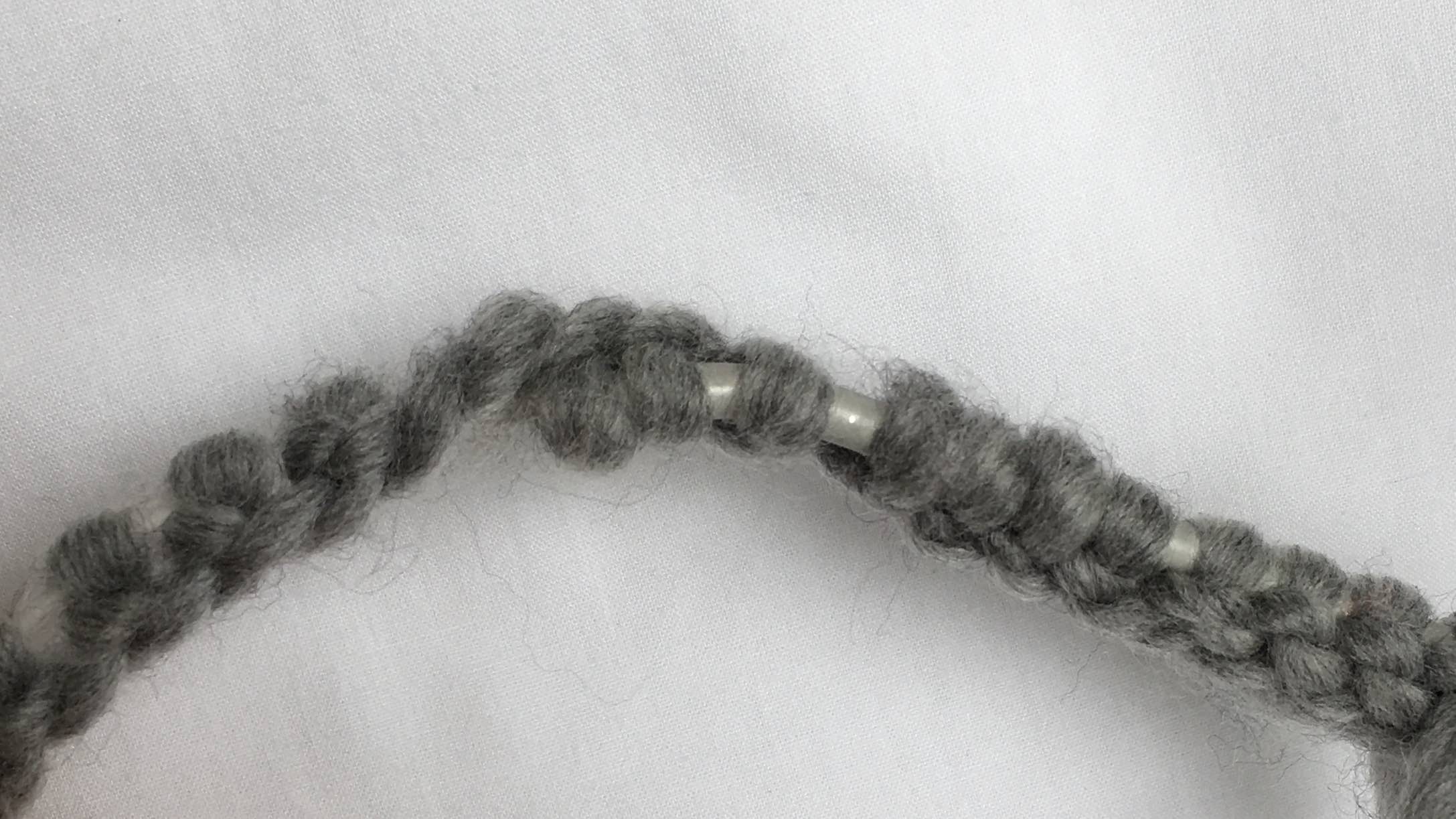}
       	    \subcaption{cast-on with $360^\circ$ twist}
        	\label{fig:1b}
    \end{minipage}
    \caption{example of circular needle with cast-on}
\label{image1}
\end{figure}

Circular knitting is a technique that uses two needles connected with a freely-moving chord. The first step of the knitting, the cast-on, runs along the chord as seen in Figure~\ref{fig:1a}. It is then closed to a circle and one then proceeds to knit a band. It is very common that the first attempt at this leads to a twisted band, similar to that seen in Figure~\ref{fig:1b}. The reason is that the cast-on has an orientation, the side with the loops and the edge. It should therefore be thought of as a ribbon of very small width. In closing the cast-on to a circle one inevitably makes a choice of the number of twists, as when gluing two ends of a ribbon. This twisting is normally viewed as a knitting fault, but we want to relate it to the notion of the framing anomaly, where one is forced to adjoin a normal vector to a closed contour and there isn't a unique choice to do that. In knitting, the normal vector is the direction in which subsequent rows are added.

\begin{figure}[h]
    \centering
    \includegraphics[width=0.42\textwidth]{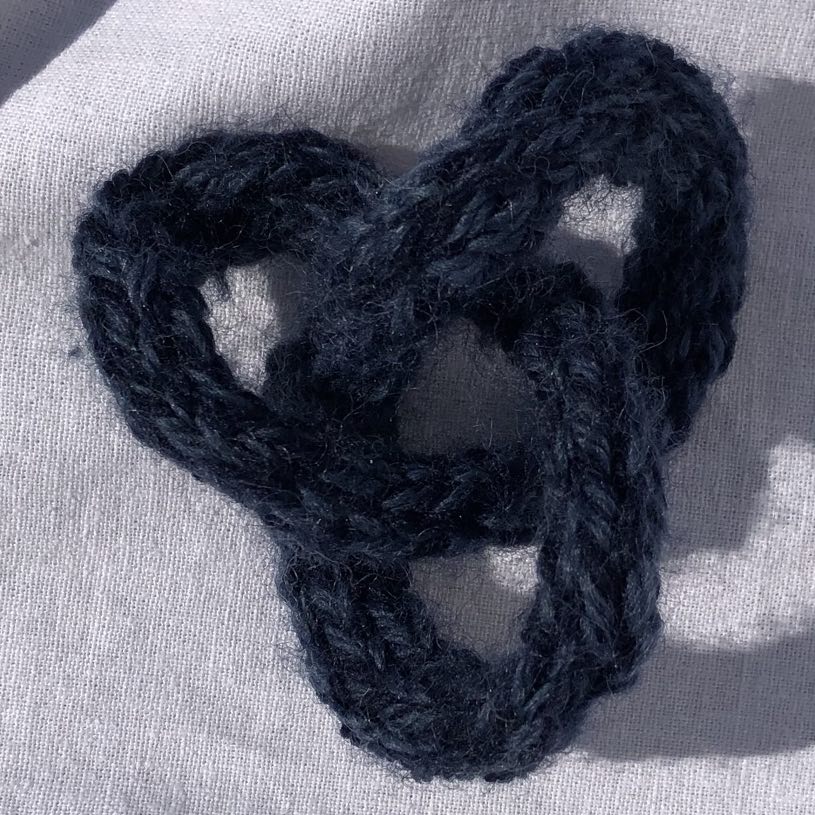}
    \caption{trefoil knot knitted using circular knitting}
    \label{fig:2}
\end{figure}

Aside from the above discussion on framing, let us comment on creating an arbitrary 
knot using circular knitting. This can be done before closing the cast-on to a circle, 
by arbitrarily knotting the chord connecting the needles. 
As this adds complication to the knitting and the maths, in the 
following we choose to focus only on the unknot (the simple circle) to highlight 
the concept of the framing alone. However, as an illustration of this process, 
we created a trefoil knot seen in Figure~\ref{fig:2}. For examples of other 
knitted knots, see \cite{Belcastro,2019_Bridges_Gallery}.

Returning to our knitted unknot. After knitting a few rows we have a band with two edges which we label $C$ (the cast-on) and $C'$ (the last knitted circle). The band itself can be seen to represent a vector field translating $C$ into $C'$. The framing number of the vector field is 
defined as the linking number \cite{Roberts_Knotes} between $C$ and $C'$. In this way we gave a mathematical definition to the twisting fault as the linking number between the two edges of the ribbon or the framing of the vector field. We illustrate this in several examples below.

The linking number is a topological invariant \cite{gauss}, and so is the framing number. Once the cast-on is closed into a a knot, the extent of twisting, or framing, is determined and can no longer be altered.

For an example, we refer back to Figure~\ref{fig:1b}, where the cast-on has a $360^\circ$ turn and is closed. Given this twisting fault, the resulting band has edges with linking number one. So our cast-on has framing one. Moreover, the arbitrary positioning of the $360^\circ$ turn along the chord results in the same framing number, which illustrates its invariant property.

Note that with regular knitting, the twist can only be an integer multiple of $360^\circ$, so one cannot create a m\"obius strip that requires a $180^\circ$ twist. There are still ways to achieve this, see for example \cite{Belcastro}, but we restrict ourselves to traditional oriented knitting.

Let us briefly mention the relation to Chern-Simons gauge theory, and for simplicity focus on 
the abelian theory with gauge group $U(1)$ \cite{witten,Grabovsky_Knotshell}. 
The Chern-Simons action on $\mathbb{R}^3$ with gauge potential $A_i$ is
\begin{equation}
    S[A]=\frac{k}{4\pi}\int_{\mathbb{R}^3}e^{ijk}A_i\partial_j A_k\,d^3x\,, \qquad k\in\mathbb{Z}\,.
\end{equation}
For a closed curve $C$, which is our knot, the Wilson loop observable is then
\begin{equation}
    W(C)= \exp\left(i\oint_{C}A_jdx^j\right).
\end{equation}
In quantum field theory one defines the expectation value of $W(C)$ as the path integral with respect to the Gaussian measure $\exp(iS[A])$. Without delving too much into the theory, Witten shows \cite{witten} that evaluating this expectation value using Feynman diagrams gives the following integral. 
\begin{equation}
   \left\langle W(C)\right\rangle = \exp\left(\frac{i}{2k}\int_{C}dx^i\int_{C}dy^j\epsilon_{ijk}\frac{(x-y)^k}{|x-y|^3}\right).
\end{equation}
One can observe that this double integral is ill-defined near $x=y$. To fix this, we need to define a copy of $C$, say $C'$ that is moved slightly away from $C$ and then compute the double integral. Now that the expression is well-defined, we can see that it essentially contains the Gauss linking integral \cite{gauss}, calculating the linking number of the two knots
\begin{equation}
\label{Gauss_Linking_integral}
    \left\langle W(C)\right \rangle=\exp\left[\frac{2\pi i}{k}\link(C,C')\right],\qquad
    \link(C,C')=\frac{1}{4\pi}\int_{C}dx^i\int_{C'}dy^j\epsilon_{ijk}\frac{(x-y)^k}{|x-y|^3}\,.
\end{equation}
Thus the calculation of the Wilson loop in Chern-Simons theory, which is related in turn to the Jones polynomial of knot theory \cite{jones}, suffers from the same problem as in circular knitting---the need to choose a normal vector field to displace the curve $C$. This is the \emph{framing anomaly}. 

This is by no means the only place in science where framing arises. DNA are famously double-stranded 
and helical, and thus DNA molecules that form closed loops, like those in the mitrochondria, 
are framed. For this and other relations of knot theory in chemistry, see for example \cite{horner}.

\section*{The Framing of an Unknot}

The framing number can be any whole number. We restricted ourselves to three knitted examples and also provide three different approaches to calculating the Gauss linking integral (\ref{Gauss_Linking_integral}) in those three cases. The examples have framing zero, one, and negative two. Note that the negative sign for the last unknot is due to opposite rotation of the twist.

\subsection*{Framing Zero}

This framing can equivalently be called standard framing or the blackboard framing. Recalling the definition, this framing constitutes a vector field that is orthogonal to the projection of the knot on the plane, oriented perpendicular to the plane. This is also the framing that is required for most knitting applications and nonzero framing is usually viewed as a fault.

\begin{figure}[h!tbp]
    \centering
    \begin{minipage}[b]{0.49\textwidth} 
	    \includegraphics[width=\textwidth]{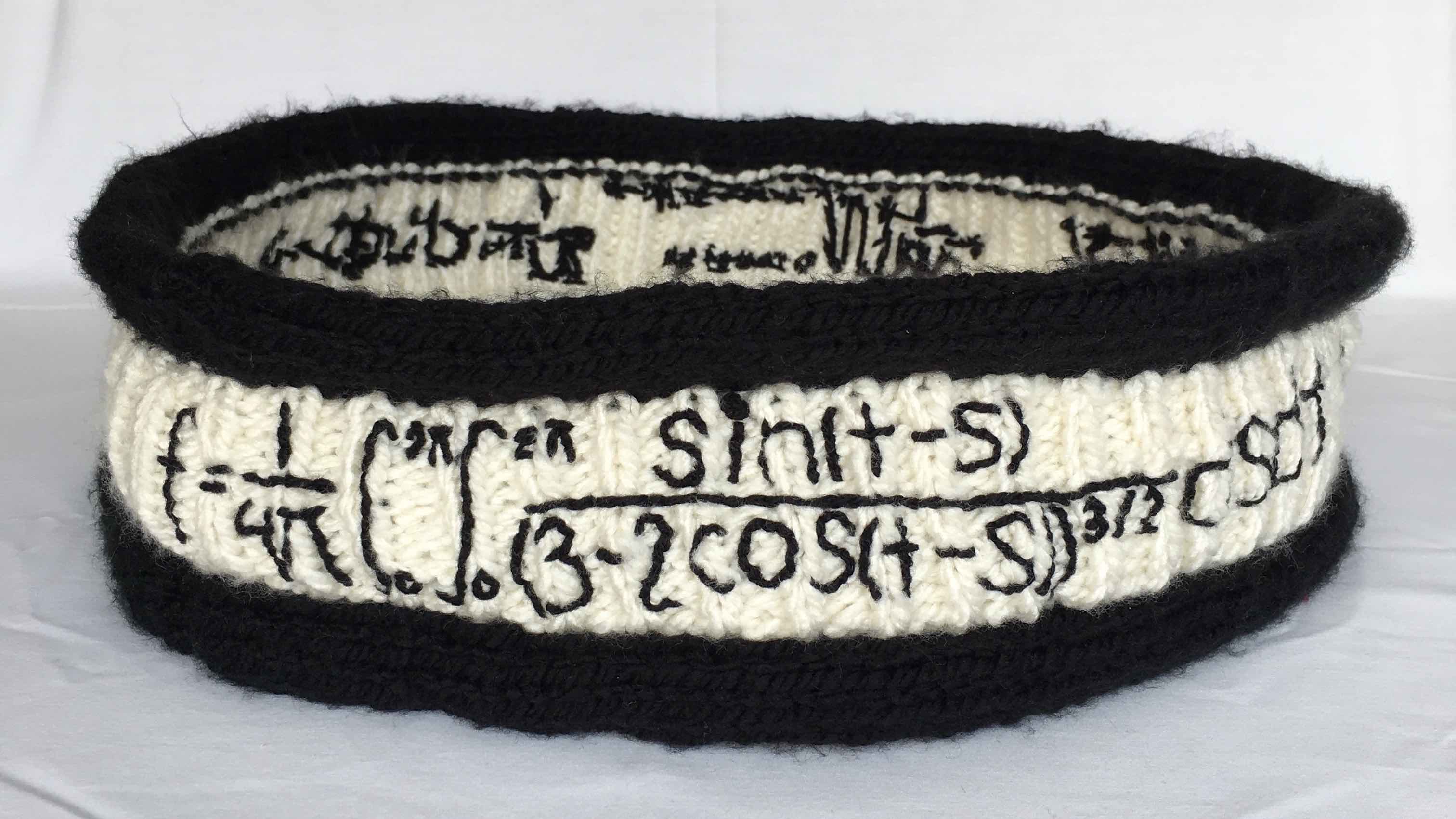}
       	    \subcaption{}
        	\label{fig:3a}
    \end{minipage}
~
    \begin{minipage}[b]{0.49\textwidth} 
	    \includegraphics[width=\textwidth]{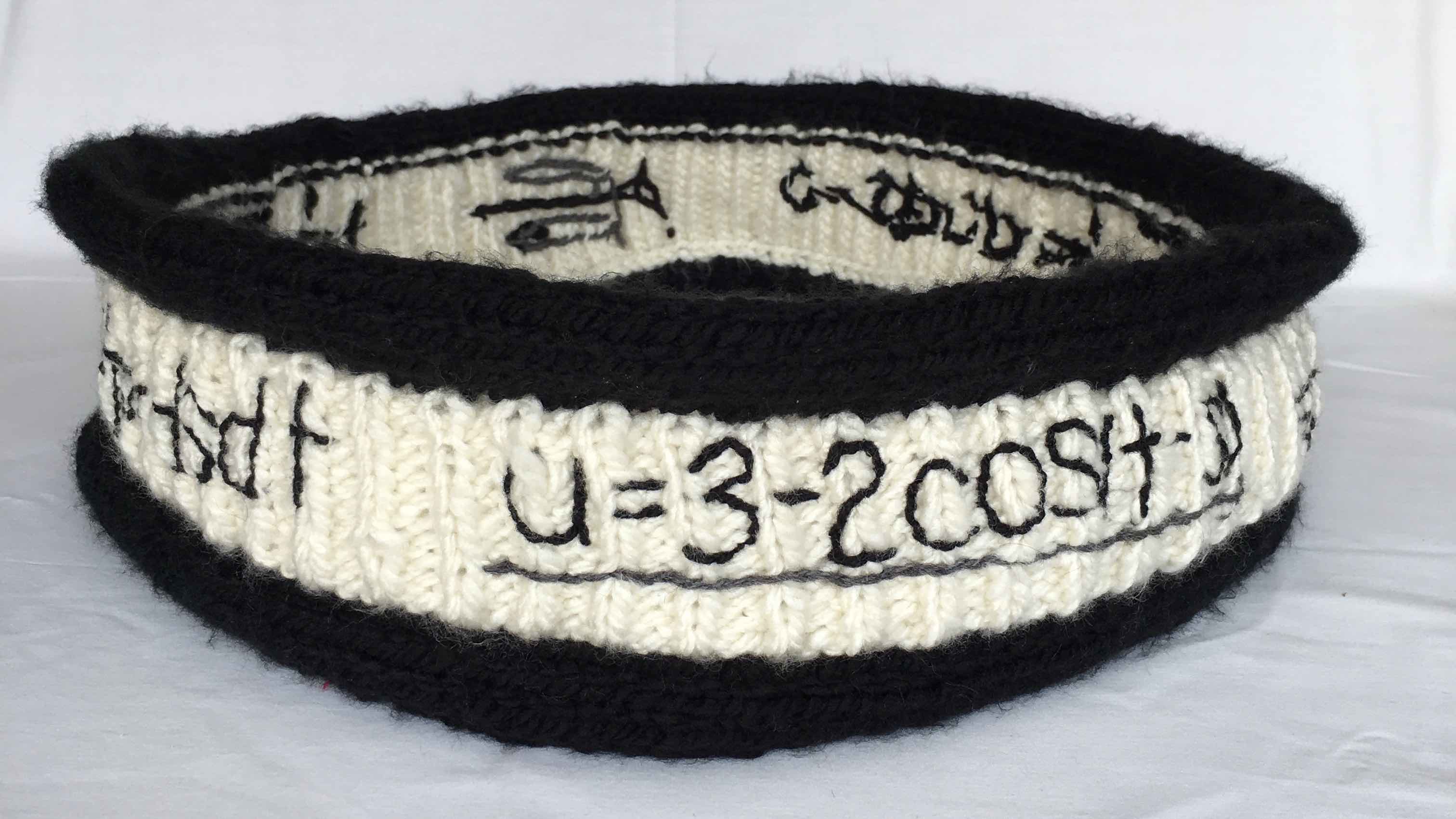}
       	    \subcaption{}
        	\label{fig:3b}
    \end{minipage}
    
    \begin{minipage}[b]{0.49\textwidth}
        \includegraphics[width=\textwidth]{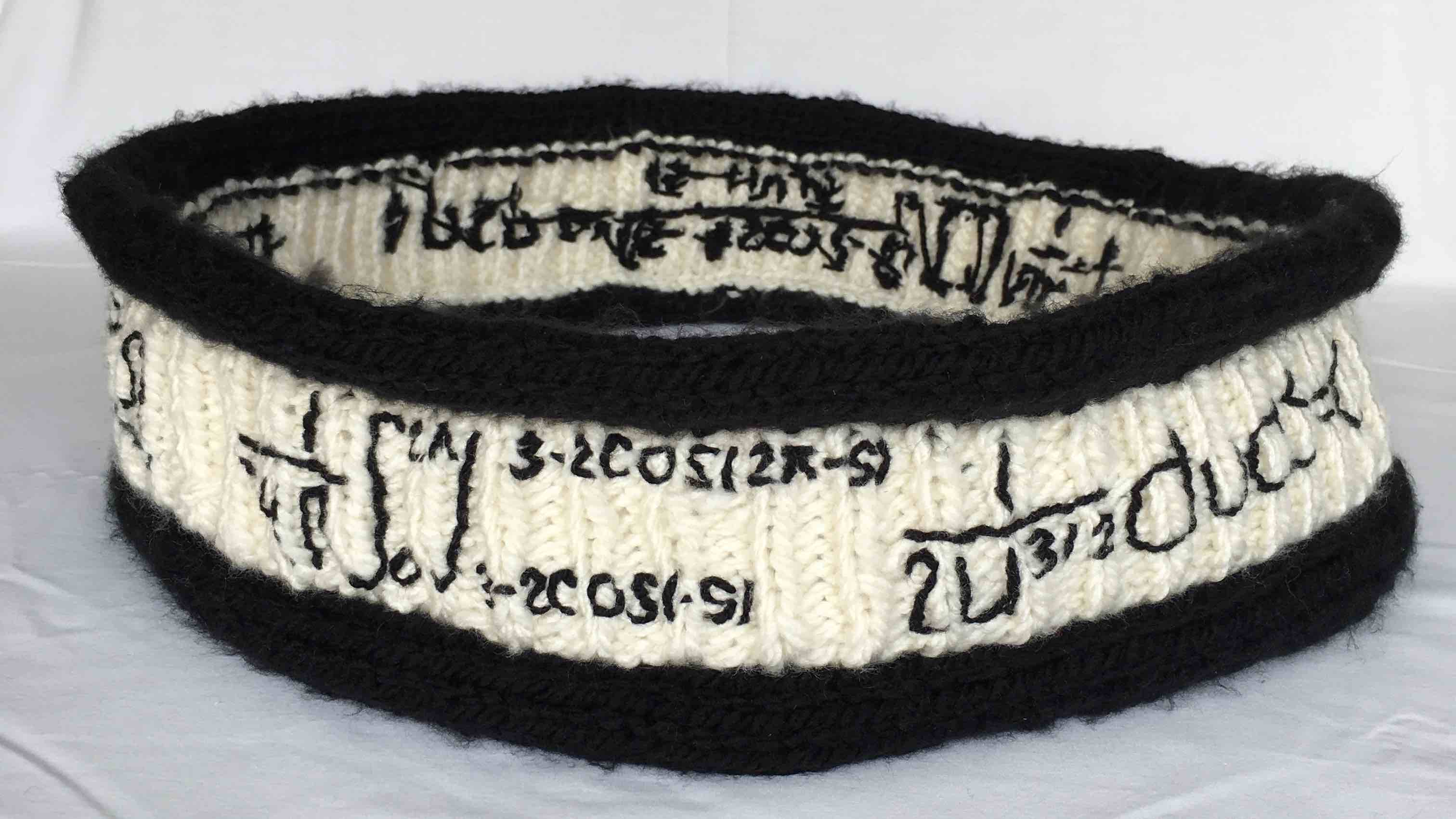}
            \subcaption{}
            \label{fig:3c}
    \end{minipage}
~
    \begin{minipage}[b]{0.49\textwidth}
        \includegraphics[width=\textwidth]{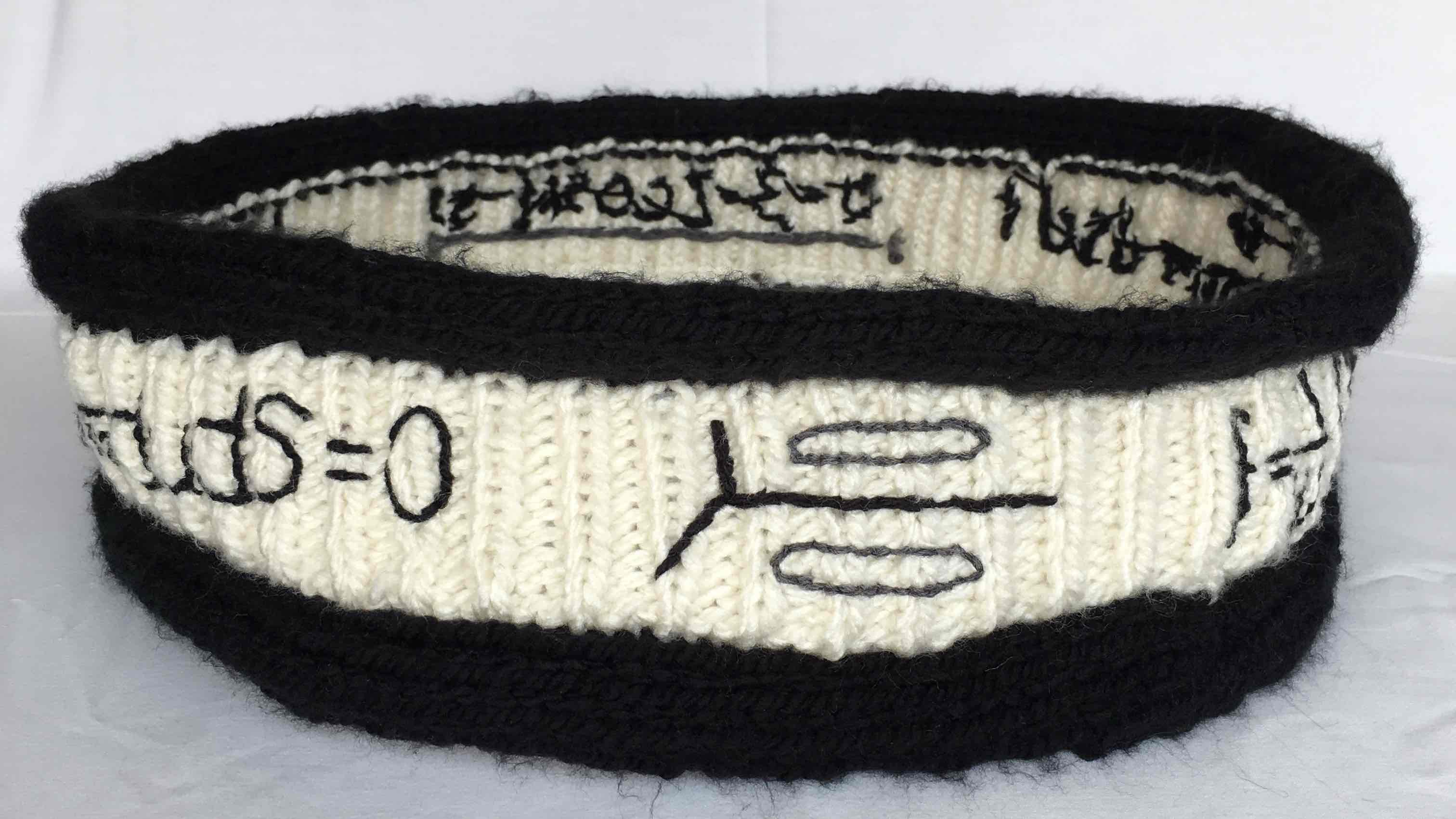}
            \subcaption{}
            \label{fig:3d}
    \end{minipage}
    \caption{$360^{\circ}$ view of unknot with framing zero}
\label{image2}
\end{figure}

The piece we created for the unknot with framing zero is shown in Figure~\ref{image2}. We provide four photographs of the unknot to show all elements of the embroidery.

For the computation of the linking integral, we need to choose a specific parametrization which we take to be
\begin{equation}
    C(t)=
    \begin{pmatrix}
        \cos(t) \\
        \sin(t) \\
        0
    \end{pmatrix},
    \qquad C'(t)=
    \begin{pmatrix}
        \cos(t) \\
        \sin(t) \\
        1
    \end{pmatrix},
    \qquad t\in [0,2\pi]\,.
\end{equation}
A graph of these curves is also embroidered on our knitted piece, see Figure~\ref{fig:3d}.

Evaluating the integrand in (\ref{Gauss_Linking_integral}), we get the double integral also seen in Figure~\ref{fig:3a}
\begin{equation}
    \link(\gamma_1,\gamma_2)=\frac{1}{4\pi}\int_{0}^{2\pi}\int_0^{2\pi}\frac{\sin(t-s)}{(3-2\cos(t-s))^{3/2}} dt ds
\end{equation}
The computation is then trivial, and uses a substitution, seen in Figure~\ref{fig:3b}. This results in the inner integral having the same boundary conditions, seen in Figure~\ref{fig:3c}, thus giving the number zero, as expected.

\noindent

Regarding the knitting, this piece is knit by starting with a variation of an I-cord cast-on 
method \cite{Suzanne_2020} in black, continued with a ribbed knit band pattern 
\cite{Slip_Stitch_Rib} using 5.5 mm circular needles and chunky acrylic yarn and finally cast-off 
with the same I-cord technique. Through the hollow black cast-on/off pieces flexible wire was threaded to increase stability and allow for shaping and posing. Embroidery is done by hand in the backstitch method, where one brings up the needle to where one wants to direct the stitch and then brings it back down to meet the previous stitch. The embroidery itself yields structural support for the band. Final dimensions are $38.1\times38.1\times12.2$ cm. 

The I chord cast-on and cast-off in black highlight the two distinct unknots. The center band connecting the two I-cords, knit in cream, is the vector field that corresponds to framing zero. The ribbed knit stitch used in this band \cite{Slip_Stitch_Rib} is chosen so as to subtly reflect the directionality of this vector field, essentially leading the viewer from the original unknot (below) to the copy (top).

\subsection*{Framing One}
\begin{figure}[h]
    \centering
    \includegraphics[width=\textwidth]{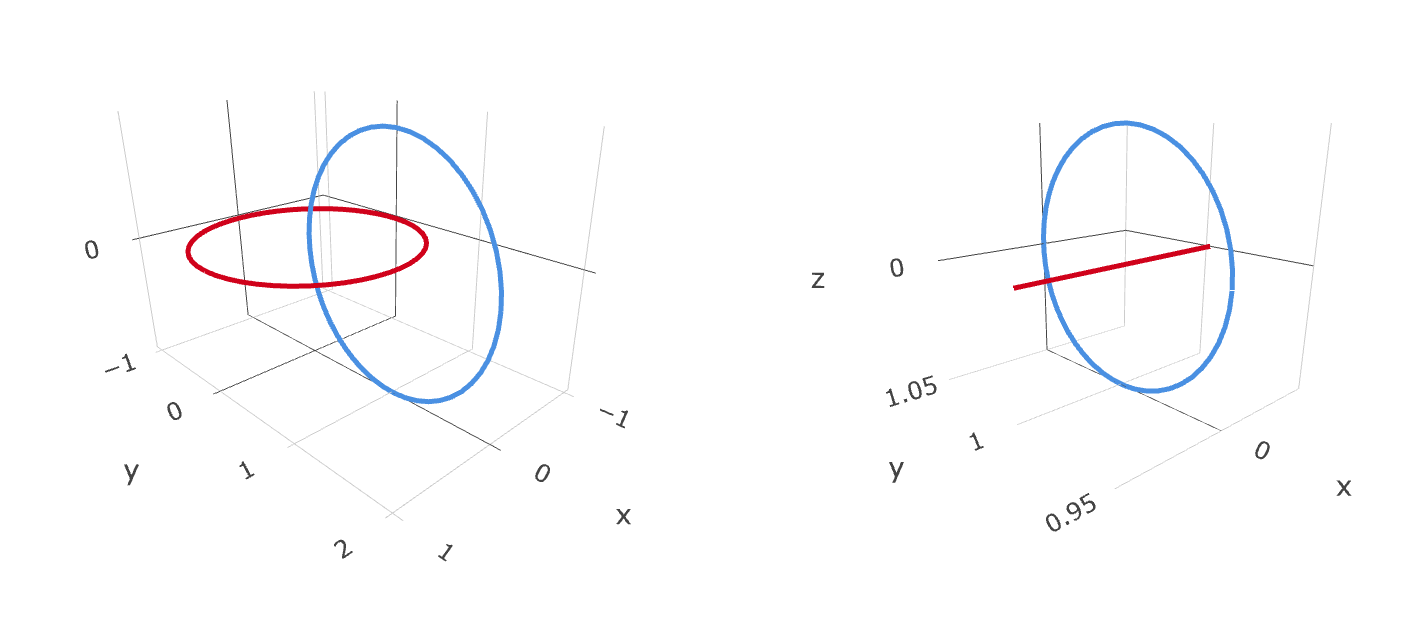}
    \caption{epsilon parametrization with $\epsilon=1$ (left) and $\epsilon=0.05$ (right)}
    \label{epsilon_fig}
\end{figure}

\begin{figure}[h!tbp]
    \centering
    \begin{minipage}[b]{0.49\textwidth} 
	    \includegraphics[width=\textwidth]{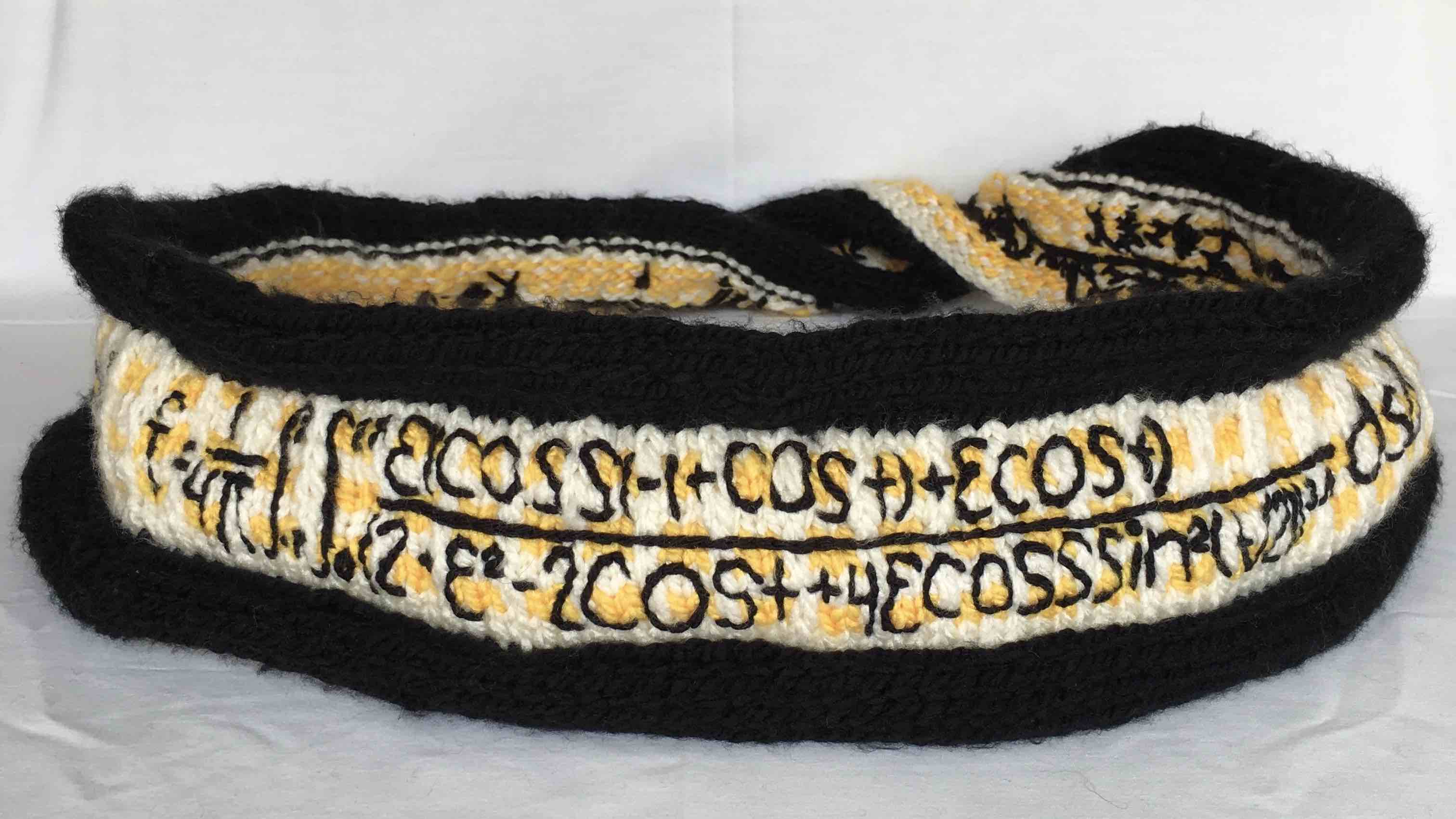}
       	    \subcaption{}
        	\label{fig:4a}
    \end{minipage}
~
    \begin{minipage}[b]{0.49\textwidth} 
	    \includegraphics[width=\textwidth]{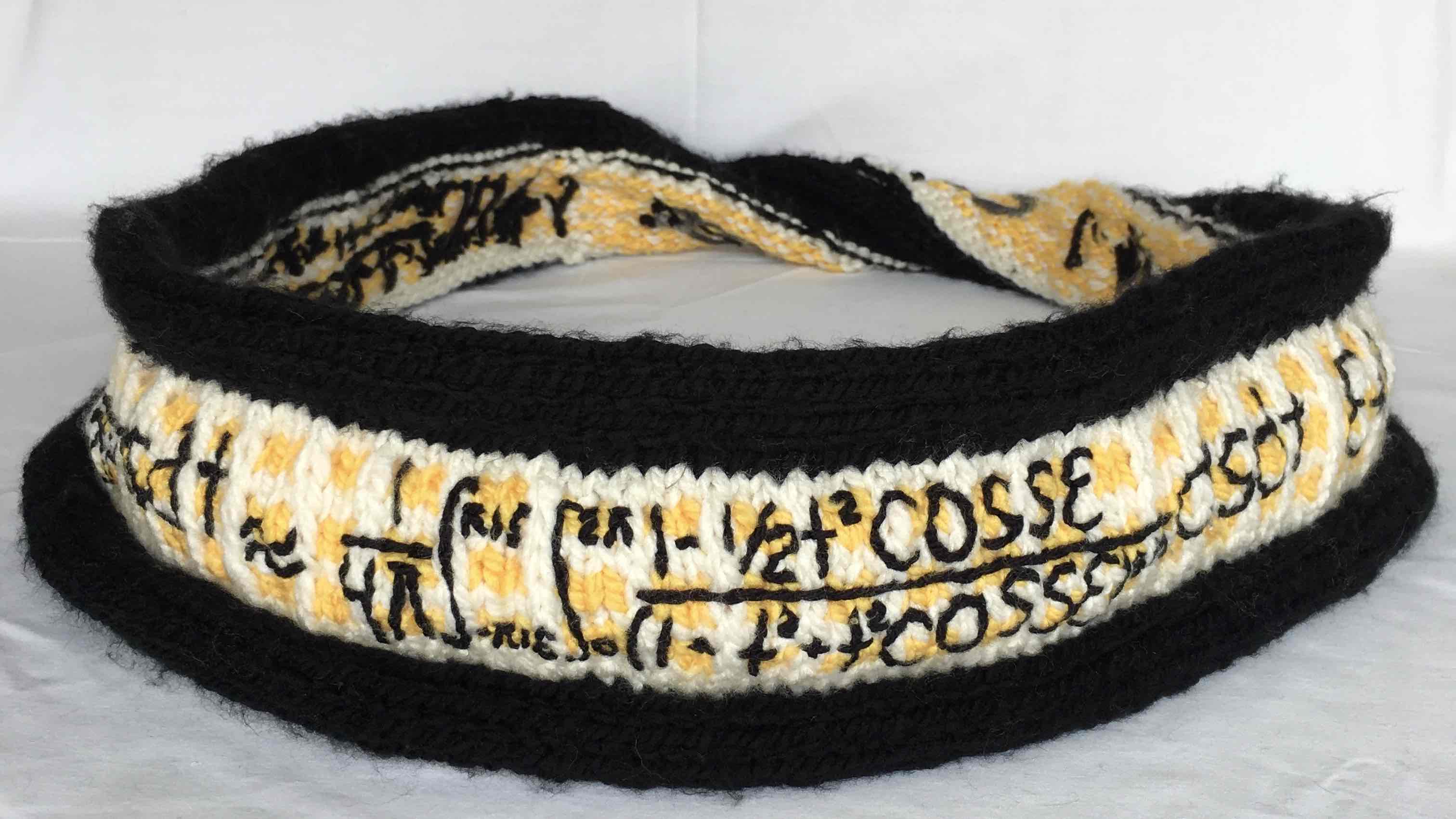}
       	    \subcaption{}
        	\label{fig:4b}
    \end{minipage}

    \begin{minipage}[b]{0.49\textwidth}
        \includegraphics[width=\textwidth]{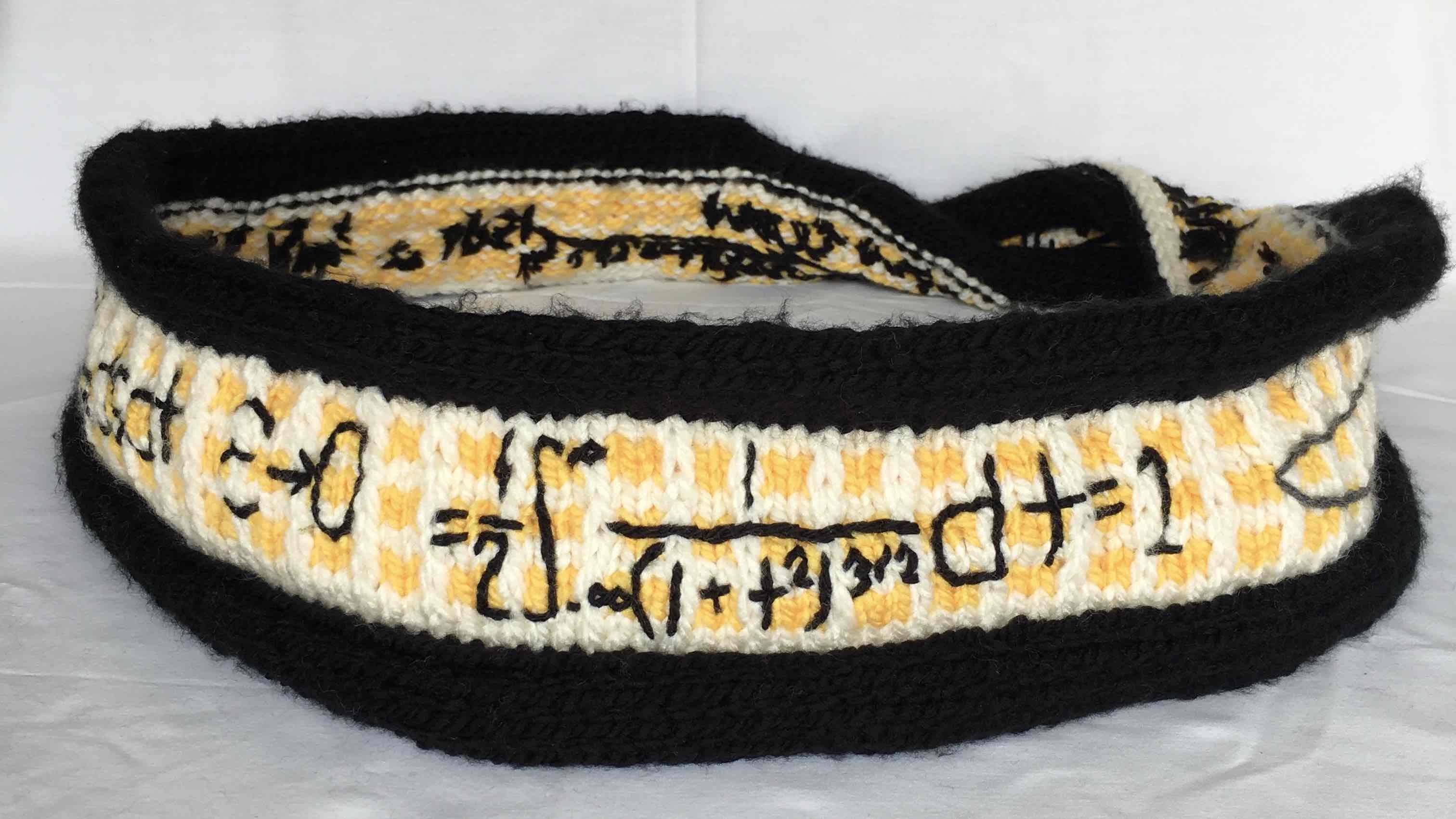}
            \subcaption{}
            \label{fig:4c}
    \end{minipage}
~
    \begin{minipage}[b]{0.49\textwidth}
        \includegraphics[width=\textwidth]{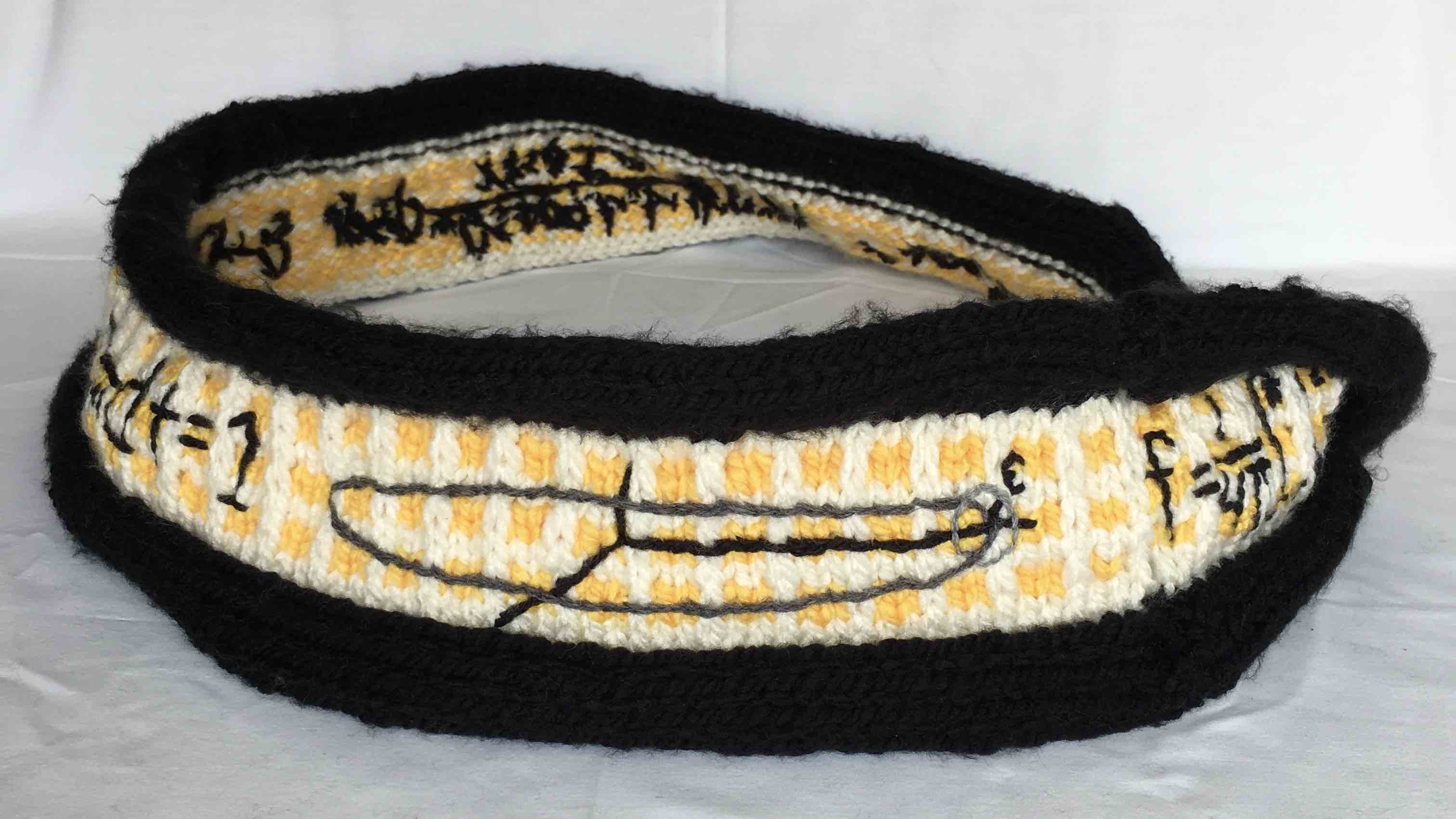}
            \subcaption{}
            \label{fig:4d}
    \end{minipage}
    \caption{$360^{\circ}$ view of unknot with framing one}
\label{image3}
\end{figure}

The knitted unknot with framing one is presented in Figure~\ref{image3}. The unknot with this framing can also be called the Hopf-link \cite{hopf}. Again, there are numerous ways to parametrize the Hopf-link. Doing the calculation for two equal size circles results in elliptic integrals, making the analytic calculation difficult. We choose instead to have one circle much smaller than the other by introducing an $\epsilon$-parameter
\begin{equation}
    C(t)=
    \begin{pmatrix}
        \sin(\epsilon t) \\
        \cos(\epsilon t) \\
        0
    \end{pmatrix},
    \qquad C'(s)=
    \begin{pmatrix}
        0 \\
        \epsilon\cos(s)+1 \\
        \epsilon\sin(s)
    \end{pmatrix},
    \qquad s\in[-\pi,\pi]\,,\quad\text{and}\quad  t\in\left[-\frac{\pi}{\epsilon},\frac{\pi}{\epsilon}\right].
\end{equation}
A graph of this is shown for two values of $\epsilon$ in Figure~\ref{epsilon_fig} and one is also embroidered in Figure~\ref{fig:4d}. As the linking number is a topological invariant, it does not depend on $\epsilon$.

Taking $\epsilon\rightarrow 0$, the circle $C'$ becomes infinitesimally small around the point $(0,1,0)$, seen as blue circle on Figure~\ref{epsilon_fig}, while the $C$ (red) circle remains of radius one, but in the vicinity of the point $(0,1,0)$ can be approximated to a line $(t,1,0)$ for $t\in(-\infty,\infty)$. 

The expression for the Gauss linking integral, also embroidered in Figure~\ref{fig:4a} is
\begin{equation}
    \link(C,C')=\frac{1}{4\pi}\int_{-\frac{\pi}{\epsilon}}^{\frac{\pi}{\epsilon}}\int_{-\pi}^{\pi} \frac{\epsilon^2(\cos(s)(\cos(\epsilon t)-1)+\epsilon\cos(\epsilon t))}{\sqrt{2+\epsilon^2-2\epsilon\cos(\epsilon t)+4\epsilon\cos(s)\sin^2\left(\frac{\epsilon t}{2}\right)}^3} ds dt
\end{equation}
We now Taylor expand the integrand for small $\epsilon$. We separately expand the numerator and the inner expression under the root of the denominator. Hence, we get (see also Figure~\ref{fig:4b})
\begin{equation}
       \frac{1}{4\pi}\int_{-\frac{\pi}{\epsilon}}^{\frac{\pi}{\epsilon}}\int_{-\pi}^{\pi}\frac{(1-\frac{1}{2}(t^2\cos(s))\epsilon+\mathcal{O}(\epsilon^2)) ds dt}{\left(1+t^2+t^2\cos(s)\epsilon+\mathcal{O}(\epsilon^2)\right)^{3/2}}
\label{equation_10}
\end{equation}
Finally, we can limit $\epsilon$ to zero to obtain 
\begin{equation}
    \frac{1}{4\pi}\int_{-\infty}^{\infty}\int_{-\pi}^{\pi}\frac{1}{(1+t^2)^{3/2}} ds dt
    =1
\label{equation_11}
\end{equation}
This is also embroidered in Figure~\ref{fig:4c}. Note, that this is the same as the Gauss linking integral (\ref{Gauss_Linking_integral}) of a line and a linked circle, seen on the right of Figure~\ref{epsilon_fig}. The calculation now is simple and we get the linking number one, corresponding to framing one.
\noindent

The same I-cord cast-on and subsequent circular knitting techniques are used as for the unknot with framing zero. Especially for the unknots with framing zero and one, the flexible wire inserted allows for easy positioning of the twists and posing of the knot to focus on different sections of embroidery. The dimensions are $38.1\times38.1\times12.9$ cm.

Furthermore, the band is knit with a bi-coloured stitch \cite{How_to_Knit} which yields more visual complexity to the piece than the framing zero one. This stitch still has the related element of the ribbing, a subtle reminder of the normal vector field.

\subsection*{Framing Negative Two}

The third and final unknot in this series is shown below in Figure~\ref{image4}. For the calculation of the framing, we choose to take two unknots that are polygons. The specific parametrization was chosen to be symmetric, which simplifies the consequent integral. We again graphed and embroidered the parametrizations, see Figure~\ref{fig:5d}. 
\begin{equation}
    C(t)=\begin{cases}
      (-1,t-2,\frac{t}{2}-1) & t\in [0,4] \\
      (t-5,2,-t+5) & t\in[4,6] \\
      (1,-t+8,\frac{t}{2}-4) & t\in[6,10] \\
      (11-t,-2,11-t) & t\in[10,12] \\
    \end{cases}
    \quad
   C'(t)=\begin{cases}
      (2-t,1,1-\frac{t}{2}) & t\in [0,4] \\
      (-2,5-t,-5+t) & t\in[4,6] \\
      (-8+t,-1,4-\frac{t}{2}) & t\in[6,10] \\
      (2,-11+t,-11+t) & t\in[10,12] \\
    \end{cases}
\end{equation}
As one can see, the integral is comprised of 16 parts, pairing segments of each path. There are relations among the 16 double integrals reducing to four expressions embroidered on the remaining parts of our knitted piece seen in Figures \ref{fig:5a}--\ref{fig:5c}. They are,
\begin{align}
\label{16_integrals}
    & \frac{1}{4\pi}\left(2\int_0^4\int_0^4\frac{2{ ds dt}}{ ((s-3)^2+(1-t)^2+(\frac{s-t}{2})^2)^{3/2}}+8\int_4^6\int_4^6 \frac{-2{ ds dt}}{((t-3)^2+(s-7)^2+(t-s)^2)^{3/2}}\right.\nonumber\\
    &+\left. 3\int_4^6\int_0^4\frac{-2{ ds ts}}{(9+(s+t-7)^2+(t+\frac{s}{2}-6)^2)^{3/2}}+3\int_0^4\int_4^6\frac{2{ ds dt}}{(1+(s-t-3)^2+(s+\frac{t}{2}-6)^2)^{3/2}}
    \right)
\end{align}

To illustrate how these integrals are evaluated it's simpler to consider two non-crossing perpendicular lines
\begin{equation}
    \gamma_1(t)=(t,0,0)\,,
    \qquad \gamma_2(t)=(0,t,1)\,.
\end{equation}
In this case we find that the integrals that take the form
\begin{equation}
    \int\int \frac{-1}{(1+s^2+t^2)^{3/2}}{ ds dt}=-\arctan\left(\frac{st}{\sqrt{1+s^2+t^2}}\right)+A(s)+B(t)\,.
\end{equation}
Similar expressions arise in our case and 
into this expression one should insert the boundary values of 
$s$ and $t$ and sum over the different connected segments. For a closed contour, 
the contributions of $A(s)$ and $B(t)$ all cancel each-other and we are left with the 
arctan's. One should be careful to account for the multi-valuedness of the 
arctan and take the correct branch cut. Doing this 
for the sum in (\ref{16_integrals}) gives $-8\pi$ and after dividing by $2\pi$, 
we get our linking number $-2$. 

\begin{figure}[H]
    \centering
    \begin{minipage}[b]{0.49\textwidth} 
	    \includegraphics[width=\textwidth]{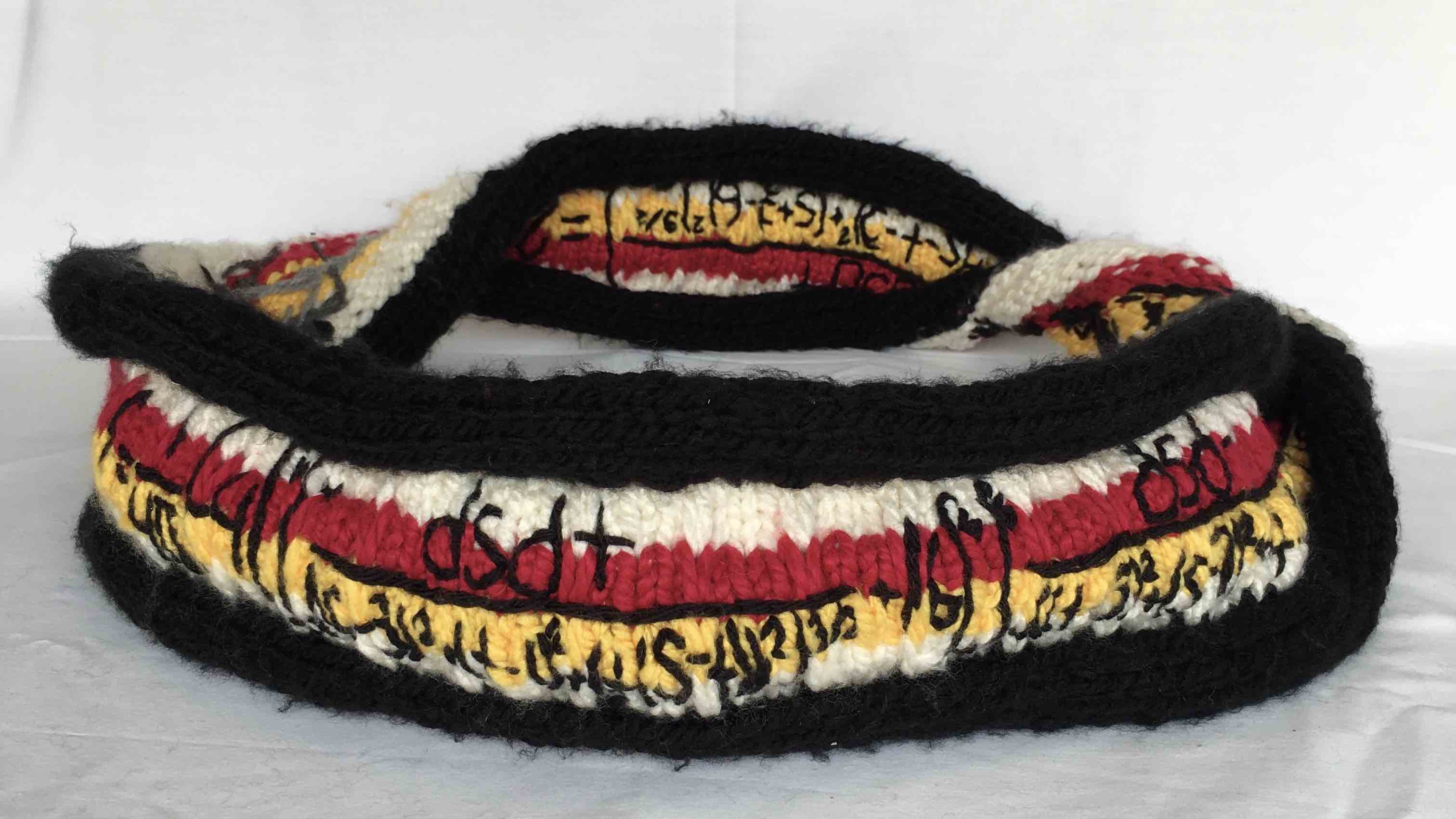}
       	    \subcaption{}
        	\label{fig:5a}
    \end{minipage}
~
    \begin{minipage}[b]{0.49\textwidth} 
	    \includegraphics[width=\textwidth]{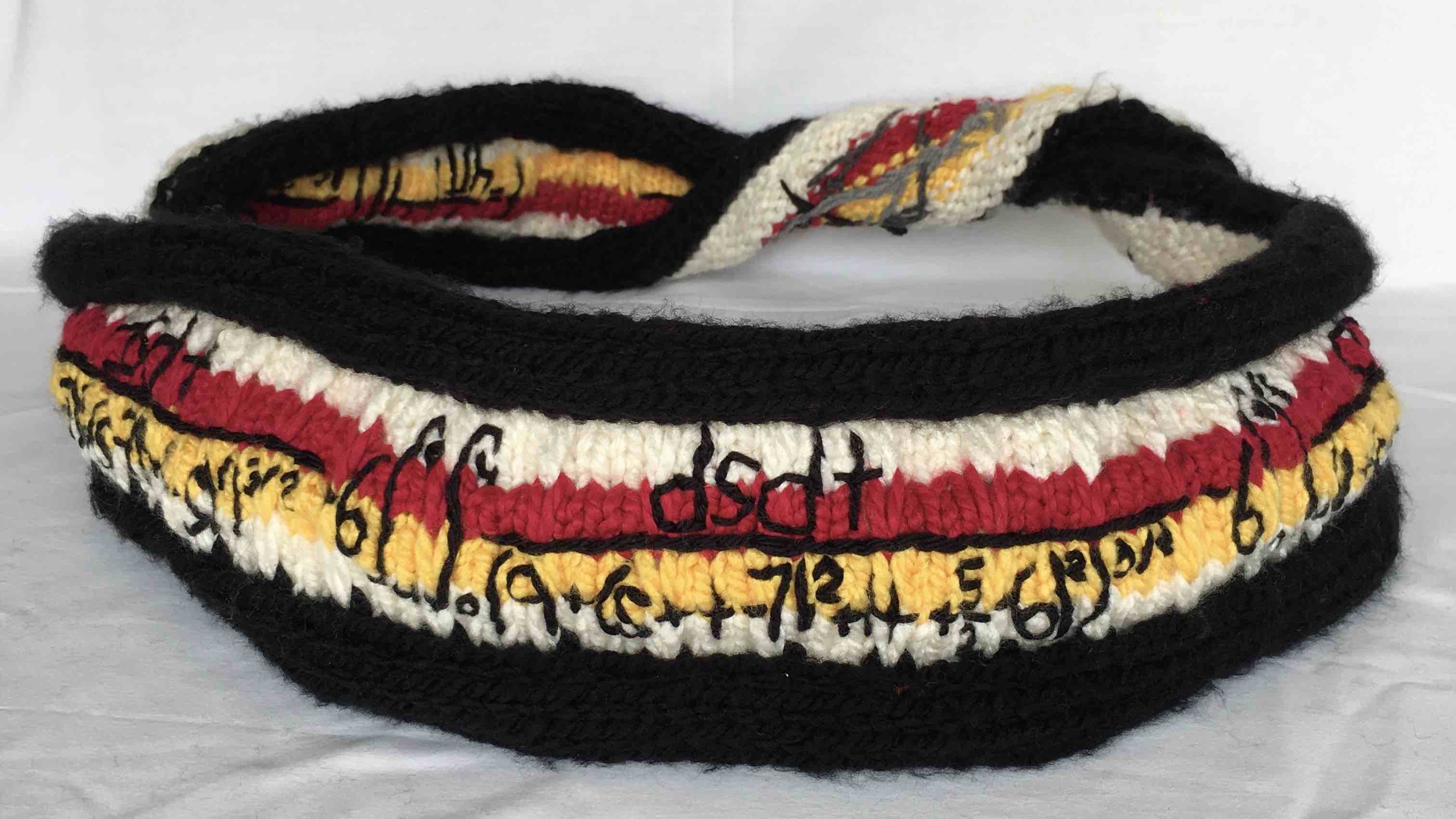}
       	    \subcaption{}
        	\label{fig:5b}
    \end{minipage}
    
    \begin{minipage}[b]{0.49\textwidth}
        \includegraphics[width=\textwidth]{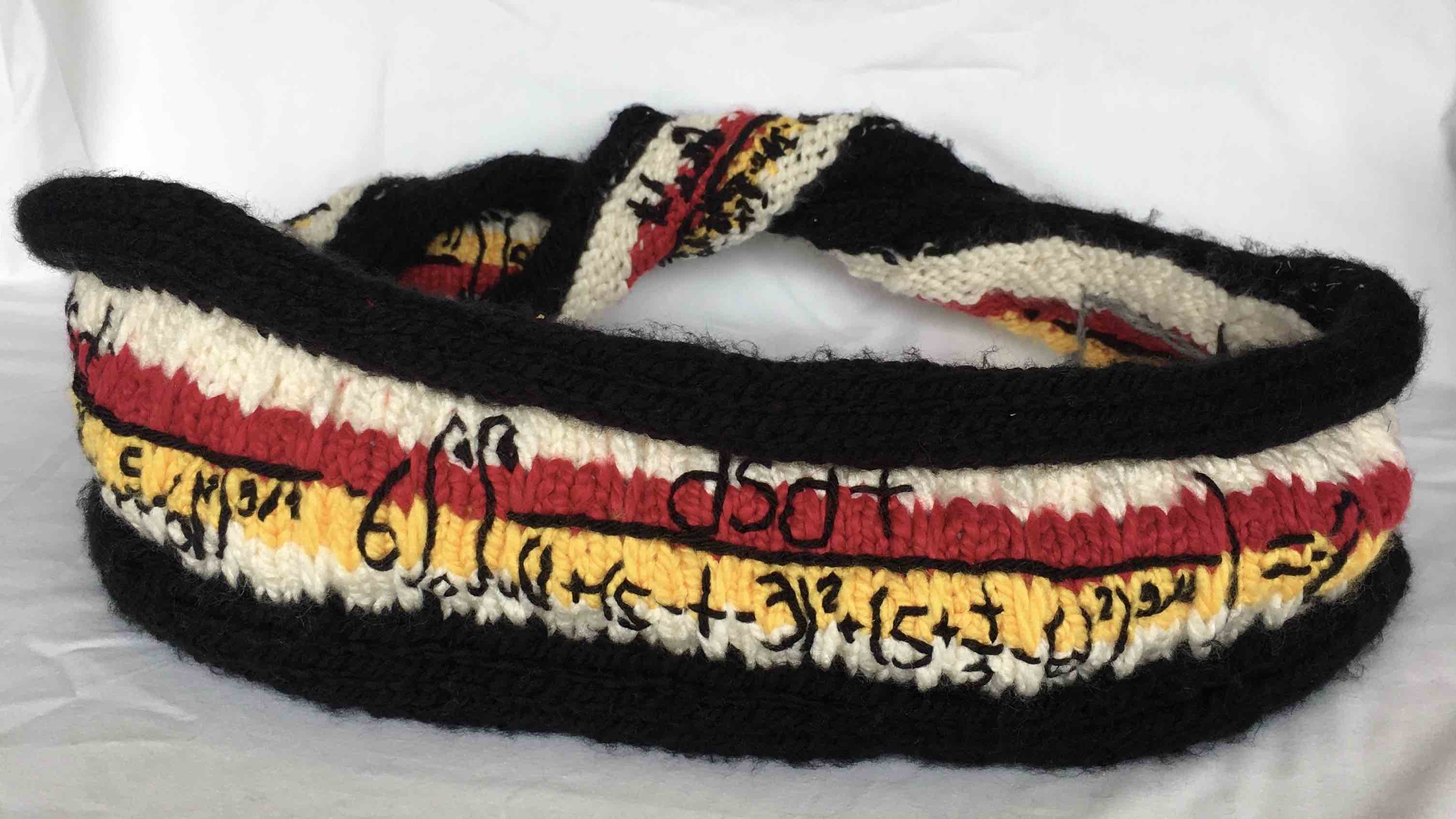}
            \subcaption{}
            \label{fig:5c}
    \end{minipage}
~
    \begin{minipage}[b]{0.49\textwidth}
        \includegraphics[width=\textwidth]{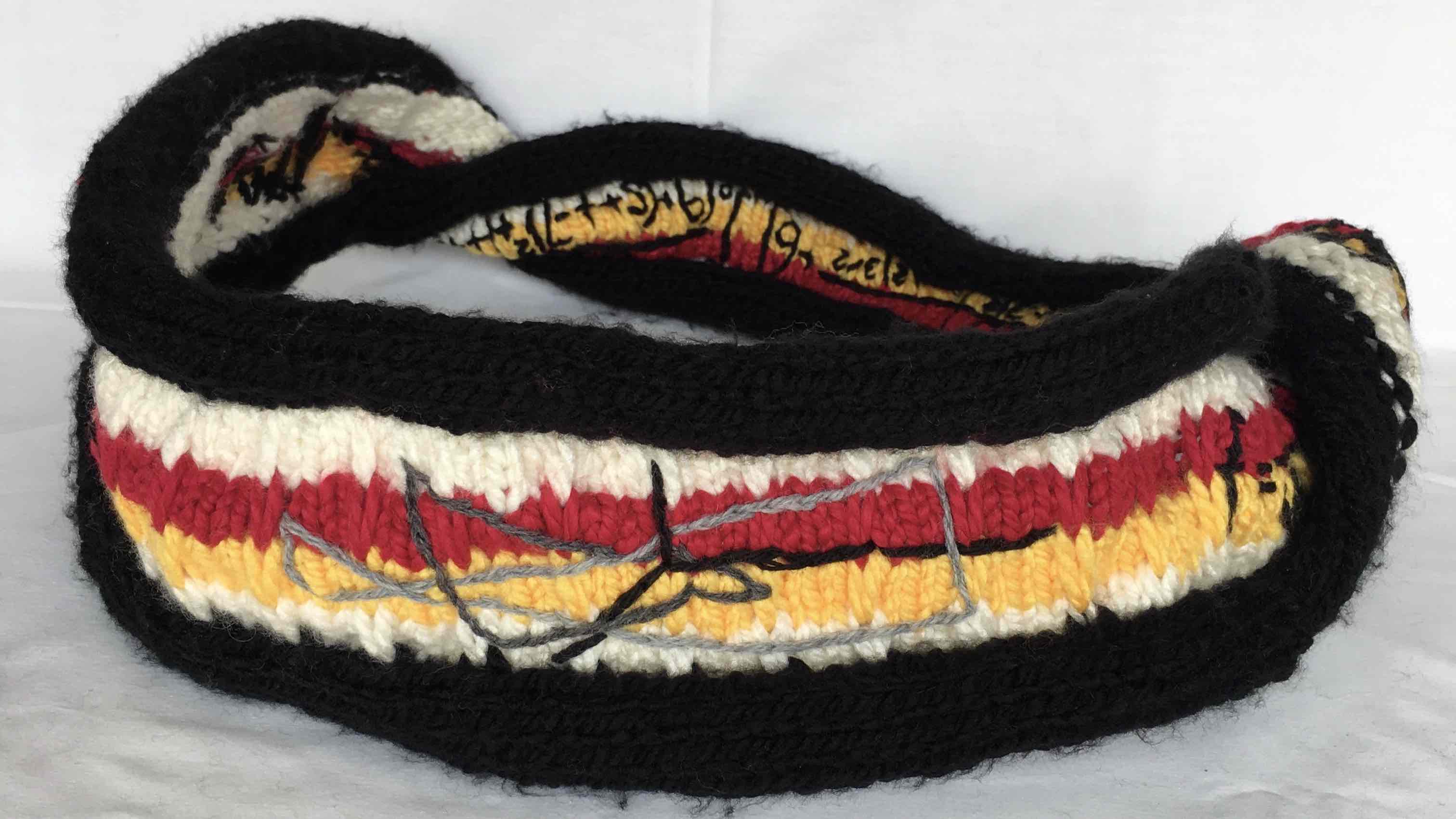}
            \subcaption{}
            \label{fig:5d}
    \end{minipage}
    \caption{$360^{\circ}$ view of unknot with framing negative two}
\label{image4}
\end{figure}

In terms of the knitting process, we use the same technique as in the previous pieces but now for the band use a tri-colored knit stitch \cite{The_Ridge_Check}. The dimensions are $38.1\times38.1\times12.5$ cm.

The amount of color in the knitting is meant to represent the increased complexity from 
framing zero with one color in the band, to one with just a slight color variation and finally to this one with framing negative 2 and the tri-color pattern.
~ 

\section*{Summary and Conclusions}
In this paper, we explore an interplay of three concepts found in knitting, mathematics 
and physics. In Witten's calculation of the Wilson loop in Chern-Simons gauge theory, 
the computation leads to the Gauss linking number which is ill-defined and requires a 
choice of framing. Likewise in circular knitting one chooses how to connect the cast-on, 
so the twisting fault can be identified with framing. 
We visualise this framing through our knitted pieces, and embroider key calculations 
that determine the framing number of a knot with given framing. 

There is always a possibility to extend this research by computing and kniting more complex knots with more complex framings, however, as we already discussed the idea remains the same. An alternative study can be done on another part of \cite{witten}, specifically the knot invariant called Jones polynomial. Again, one can attempt to knit some knots and embroider the calculation of Jones polynomial on the knot itself. 

\section*{Acknowledgements}
First and foremost we would like to thank King's Undergraduate Research Fellowship (KURF) team for facilitating and funding this research project. Moreover, we would like to thank Angel Paznokas for her advice, support, and her knitting expertise throughout this project. 

    
{\setlength{\baselineskip}{13pt} 
\raggedright				

} 
   
\end{document}